\input amssym.def
\input amssym.tex
 
\input xy
\xyoption{all}

\magnification=\magstep1

\hfuzz=10pt
\font\bigrm=cmb10 scaled 1200
\font\nrm=cmcsc10 at10pt

\def\ra{\rightarrow}
\def\lra{\longrightarrow}
\def\geq{\geqslant}
\def\leq{\leqslant}
\def\fini{{$\quad\quad\square$}}
\def\g{{\gamma}}
\def\s{{\sigma}}
\def\a{{\alpha}}
\def\de{{\delta}}
\def\dd{{\delta}}
\def\im{{\frak m}}
\def\ib{{\frak b}}
\def\ip{{\frak p}}
\def\proj{{\rm Proj}}
\def\codim{{\rm codim}}

\def\reg{{\rm reg}}
\def\ann{{\rm Ann}}

\def\tor{{\rm Tor}}
\def\jac{{\rm Jac}}
\def\Z{{\cal Z}}
\def\Sc{{\cal S}}
\def\zer{{\cal Z}}
\def\azer{{\frak K}}
\def\acc{{\widehat{\frak C}}}
\def\ac{{\frak C}}
\def\azerc{{\widehat{\frak K}}}

\def\breg{\ib\hbox{\rm -reg}}

\centerline{\bigrm Cohomology of  projective schemes: From annihilators
  to vanishing}
\bigskip
\centerline{Marc Chardin}

\bigskip\bigskip\medskip
\centerline{\bf Introduction}\bigskip

This article comes from our quest for bounds on the
Castelnuovo-Mumford 
regularity of schemes in terms of their ``defining equations'', in the
spirit of [BM], [BEL], [CP] or [CU]. The references [BS], [BM], [V] or [C] 
explains how this notion of regularity is a mesure of the algebraic
complexity of the scheme, and provides several computational
motivations.

It was already remarked by several authors (see for instance [M], [MV],
[NS1] or [NS2]) that one may bound the Castelnuovo-Mumford regularity
of a Cohen-Macaulay projective scheme in terms of its $a$-invariant and
the power of the maximal ideal that kills all but the top local
cohomology modules. Such a connection is a particular case of Lemma
2.0, which shows the way but will not be used in the sequel, because
some natural killers have stronger annihilation properties that leads
to sharper estimates.   

In connection with our previous joint work with Philippon, we
introduce partial annihilators of modules ({\it i.e.} elements that
annihilates in some degrees) and prove that uniform
(partial)annihilators of Koszul homology modules give rise to
(partial)annihilators of \v Cech cohomology modules of the quotient by a
sequence of parameters (Proposition 2.1). Combined with 
Proposition 2.3, this leads to our key for passing from annihilators to
vanishing: Proposition 2.4.

The third section gathers results on uniform annihilators that
have two main sources: tight closure and liaison. For the applications
to Castelnuovo-Mumford regularity, the key is to determine an
annihilator of the cohomology modules on which we have a control in
terms of degrees of generators. The Jacobian ideal, which kills
phantom homology by a theorem of Hochster and Huneke, and the ideal we
construct via liaison are the only ones that we are able to control
today, it would be interesting to have such a control on other natural
annihilators (e. g. the parameter test ideal). 

The results on regularity are combinations of the two preceeding
ingredients: control of annihilators and passing from annihilators to
vanishing. They hold in any characteristic and, in positive 
characteristic, improves the ones of [CU] for the unmixed part
(even if not completely comparable). Also, they do not rely on Kodaira
vanishing. The main result of a new type is the following (Theorem 4.4),\medskip

{\bf Theorem.} {\sl Let $R$ be a polynomial ring over a field,
$I\subset R$ be a homogeneous $R$-ideal generated by forms of degrees
$d_{1}\geq\cdots \geq d_{s}$. Consider $r\leq s$ and  $J\subset R$,
an intersection of isolated primary component of codimension $r$ of
$I$. Set $\Sc :=\proj (R/J)$,  $\Z :=\proj (R/I)$ and assume that  $\Z
\not= \emptyset$ and

{\rm (1)}  $\Z$ have at most isolated singularities on
$\Sc$,\smallskip 

{\rm (2)} $\Sc$ do not meet the other components of $\Z$.
\smallskip
Then
$$
\reg (R/J)\leq \dim (R/J)(d_{1}+\cdots +d_{r}-r-1)+1.
$$
}

It is a first step in the direction of bounding the regularity of all
the isolated components of a scheme in terms of degrees of generators.
It seems reasonable to hope that hypothesis (2) may be removed; it may
even be that (1) and (2) are superfluous. 

Most of the estimates we are aware of gives either huge bounds or only
concerns the 
top dimensional component. There are two exceptions: our result with
Philippon on the zero dimensional part (which is also crucial here)
and our joint work with Ulrich [CU, 4.7 (b)], where we have an
hypothesis on the singularities of the total scheme $\Z$ in place of
the component $\Sc$.  
\bigskip

{\bf \S 1. Notations}\medskip
Set $R:=k[X_{0},\ldots ,X_{n}]$ and $\im :=(X_{0},\ldots ,X_{n})$, where $k$ 
is a field.\medskip

If $\g =(\g_{1},\ldots ,\g_{s})$ is a collection of elements of $R$,
and $M$ is a $R$-module
$K_{\bullet}(\g \, ;M)$ will denote the Koszul complex of $\g$ on $M$
and $H_{i}(\g \, ;M)$ its $i$-th homology module (and similar
notations for Koszul cohomology). 
We denote by ${\cal C}^{\bullet}_{\im}(M)$ the \v Cech complex
$$
0\lra M\lra \oplus_{i}M_{X_{i}}\lra
\oplus_{i<j}M_{X_{i}X_{j}}\lra\cdots \lra M_{X_{0}\cdots X_{n}}\lra 0
$$
(starting at ${\cal C}^{0}_{\im}(M)=M$) and by $H^{i}_{\im}(M)$
its $i$-th cohomology module.\medskip

For a graded $R$-modules $M$, we will set 
$a_{i}(M):=\sup \{ \mu\ |\ H^{i}_{\im}(M)_{\mu}\not= 0\}$ and $b_{i}(M):=\sup \{ \mu\ |\
\tor_{i}^{R}(M,k)_{\mu}\not= 0\}$.
If $M$ is finitely generated, $a_{i}(M)$ is finite and therefore
$$
\reg (M):=\max_{i}\{ a_{i}(M)+i\} =\max_{i}\{ b_{i}(M)-i\}
$$ 
also is. The $a$-invariant of $M$ is  $a(M):=a_{\dim M}(M)$,
and ${\rm Tor}_{i}^{R}(M,k)
=H_{i}(X_{0},\ldots ,X_{n}\, ;M)$.\medskip

If $B$ is a standard graded algebra, $H$ a graded $B$-module and
$S\subseteq {\bf Z}$, we set 
$$
\ann_{B}^{S}(H):=\{ x\in B\ |\ (xH)_{\nu}=0,\ \forall \nu\in S\}.
$$

We denote by $\zer_{M}$ the set of tuples $z$ of
homogeneous elements such that $\dim M/(z)M=0$, set $d:=\dim M$ and 
$$
{\azer}^{\mu}_{M}:=\bigcap_{z\in\zer_{M}
  }\bigcap_{i}\ann^{>\mu -i}_{R}(H^{i}(z ;M)), \quad
{\azerc}^{\mu}_{M}:=\bigcap_{z\in\zer_{M}
  }\bigcap_{i<d}\ann^{>\mu -i}_{R}(H^{i}(z ;M)),
$$
%and also
$$
{\ac}^{i,\mu}_{M}:=\ann^{>\mu -i}_{R}(H^{i}_{\im}(M)), \quad
{\ac}^{\mu}_{M}:=\bigcap_{i}{\ac}^{i,\mu}_{M}, \quad
{\acc}^{\mu}_{M}:=\bigcap_{i<d}{\ac}^{i,\mu}_{M}. 
$$

If $M$ is Cohen-Macaulay and equidimensionnal,
${\acc}^{\mu}_{M}={\azerc}^{\mu}_{M}=A$ for any $\mu$. 
With these definitions,
$$
\reg (M)=\min \{ \mu\ |\ {\ac}^{\mu}_{M}=A\}=\min \{ \mu\ |\
{\azer}^{\mu}_{M}=A\} .
$$
Following [CP] we define $\breg (M):=\min \{ \mu\ |\ \ib\subseteq
{\azer}^{\mu}_{M}\}$. \bigskip

{\bf \S 2. Partial annihilators of Koszul and local cohomologies}
\medskip
Before investigating results that are better fitted to the kind of
annihilators we are able to build in the next paragraph, let us state
a basic result which is typical of how one passes from annihilators to
vanishing,\medskip

{\bf Lemma 2.0.} {\sl Let $M$ be a finitely generated graded $R$-module of 
dimension $d$, $\g =(\g_{1},\ldots ,\g_{s})$ a collection of homogeneous 
elements of $R$ of degrees $\dd_{1}\geq \cdots \geq \dd_{s}$ and $J$ the $R$-ideal 
they generate.

Let $N,i$ be integers and assume that,\smallskip
{\rm (1)} $J\subseteq \ann^{> N}_{R}(H^{i}_{\im}(M))$,\smallskip
{\rm (2)} $a_{i+j}(H_{j}(\g \, ;M))\leq N$, for every $j$,\smallskip
{\rm (3)} $a_{i+k}(M)\leq N-(\dd_{1}+\cdots +\dd_{k+1})$ for $0<k<s$.\smallskip
Then,
$$
a_{i}(M)\leq N.
$$
} 
To apply the result, note that if, for example, $M$ is Cohen-Macaulay
on the punctured spectrum and $0<i<d$, choosing for $\g$ a complete
system of homogeneous parameters in $\cap_{i\not= d}\ann
(H^{i}_{\im}(M))$, (1) and 
(2) are void, and (3) is obtained by descending recursion
on $i$ from an estimate of the $a$-invariant of $M$.\medskip
{\bf Proof.} Consider the double complex, ${\cal C}^{\bullet}_{\im}
K_{\bullet}(\g\, ;M)$. It gives rise to two spectral sequences, whose 
first terms are ${'E}_{1}^{pq}={\cal C}_{\im}^{p}(H_{q}(\g\, ;M))$ 
and ${'E}_{2}^{pq}=H^{p}_{\im}(H_{q}(\g\, ;M))$, for the first one. 

If $\mu >N$ and $p-q=i$, $({'E}_{2}^{pq})_{\mu}=0$ by (2), and
therefore $({'E}_{\infty}^{pq})_{\mu}=0$.

On the other hand the second spectral sequence gives, 
$$
{''E}_{1}^{pq}=K_{q}(\g\, ;H^{p}_{\im}(M))=\bigoplus_{i_{1}<\cdots <i_{q}}
H^{p}_{\im}(M)[-\dd_{i_{1}}\cdots -\dd_{i_{q}}].
$$

Note that ${''E}_{1}^{i0}=H^{i}_{\im}(M)$ and 
${''E}_{2}^{i0}=H^{i}_{\im}(M)/JH^{i}_{\im}(M)$ so that $({''E}_{2}
^{i0})_{\mu}=H^{p}_{\im}(M)_{\mu}$ for $\mu >N$ by (1). 

If $k>0$, $({''E}_{1}^{i+k,k+1})_{\mu}=0$ for $\mu >N$ by (3) (also note that 
$E_{0}^{pq}=0$ for $q>s$). As a consequence, 
$({''E}_{k+1}^{i+k,k+1})_{\mu}=0$ for $k>0$ and $\mu >N$, so
that $0=({''E}_{\infty}^{i0})_{\mu}=({''E}_{2}^{i0})_{\mu}=H^{i}_{\im}(M)_{\mu}$
if $\mu >N$.\fini\medskip

We now turn to our main technical result on annihilators,\medskip

{\bf Proposition 2.1.} {\sl Let $M$ be a finitely generated graded
  $R$-module of dimension $d$. Then ${\azer}^{\mu}_{M}\subseteq
  {\ac}^{\mu}_{M}$ and   ${\azerc}^{\mu}_{M}\subseteq
  {\acc}^{\mu}_{M}$. Furthermore, if  $\g =(\g_{1},\ldots ,\g_{s})$ is
  an homogeneous system of parameters, $\varepsilon :=\dim 
  H_{1}(\g\, ;M)$ and $\s :=\sum_{i=1}^{s} (\deg \g_{i}-1)$, then
\smallskip
{\rm (i)} For $\varepsilon\leq i<d-s$,
${\azerc}^{\mu}_{M}\subseteq {\ac}^{i,\mu +\s}_{M/(\g )M}$. 
\smallskip
{\rm (ii)} For $i\geq \varepsilon$, ${\azer}^{\mu}_{M}\subseteq
  {\ac}^{i,\mu +\s}_{M/(\g )M}$.
}\medskip

For the proof, we will need the following lemma, which is part of the
``folklore'' (see e. g. [V, 8.3.1] for a similar statment),\medskip 

{\bf Lemma 2.2.} {\sl Let $B$ be a standard graded algebra, $\im :=B_{>0}$,
$x:=(x_{1},\ldots ,x_{d})$ forms in $B_{1}$ and set $x^{t}:=(x_{1}^{t},
\ldots ,x_{d}^{t})$.
For every $i$ and $t>0$, and every finitely generated graded $B$-module $M$
such that $\dim (M/(x)M)=0$, there exists a homogeneous morphism of 
$B$-modules,
$$
\psi^{i}_{t}:H^{i}(x^{t};M)\lra H^{i}_{\im}(M)
$$ 
which is functorial in $M$ and such that the restriction $\psi^{i}_{t,\mu}:
H^{i}(x^{t};M)_{\geq\mu}\lra H^{i}_{\im}(M)_{\geq\mu}$ is an isomorphism for
$t >\max_{i}\{ a_{i}(M)\} -\mu$.}
\medskip

Note that the condition  $\dim (M/(x)M)=0$ is always satisfied 
if $\dim (B/(x))=0$.\medskip

{\bf Proof of Lemma 2.2.} Consider the two spectral sequences arising from
the double complex $D^{\bullet\bullet}:={\cal C}^{\bullet}_{\im}
K^{\bullet}(x^{t};M)$. 
As  $\dim (M/(x)M)=0$, one spectral sequence abouts at the first step 
and provides an isomorphism $H^{q}(D^{\bullet}_{tot})\simeq 
{'E}_{1}^{0q}=H^{q}(x^{t};M)$.

On the other hand, the second spectral sequence have as first terms:
$$
{''E}_{1}^{pq}=H_{\im}^{p}(M)[tq]^{{q}\choose{d}}
$$
and one has inclusions ${''E}_{1}^{p0}\supseteq {''E}_{2}^{p0}\cdots 
\supseteq {''E}_{\infty}^{p0}$. We define $\psi^{i}_{t}$ 
as the composed map,
$$
\xymatrix{
H^{i}(x^{t};M)= H^{i}(D^{\bullet}_{tot})\ar@{->>}[r]& {''E}_{\infty}^{i0}
\ar@{^{(}->}[r]& {''E}_{1}^{i0}=H_{\im}^{i}(M).\\}
$$

To see that $\psi^{i}_{t,\mu}:=(\psi^{i}_{t})_{\geq\mu}$ is an isomorphism for 
$t >\max_{i}\{ a_{i}(M)\} -\mu$, remark that for such $t$ and $\mu$
one has $({''E}_{1}^{pq})_{\mu}=0$ for $q>0$, so that all the above
inclusions are equalities and $H^{i}(D^{\bullet}_{tot})_{\mu}\ra 
({''E}_{\infty}^{i0})_{\mu}$ is an isomorphism.\fini\medskip

{\bf Proof of Proposition 2.1.} We may assume that $k$ is infinite.

For the first statment, consider a sequence of parameters $y_{1},\ldots ,y_{d}$ in
$A_{1}$. By hypothesis $[ {\azer}^{\mu}_{M}H^{i}(y_{1}^{t},\ldots
,y_{d}^{t};M)]_{> \mu -i}=0$ for   
any $i$ and $t>0$; therefore  $[{\azer}^{\mu}_{M}H^{i}_{\im}
(M)]_{> \mu -i}=0$ for any $i$ by Lemma 2.2. Also $[
{\azerc}^{\mu}_{M}H^{i}(y_{1}^{t},\ldots ,y_{d}^{t};M)]_{> \mu -i}=0$ for  
any $i\not= d$ and $t>0$ implies that  $[{\azerc}^{\mu}_{M}H^{i}_{\im}
(M)]_{> \mu -i}=0$ for $i\not= d$. 

For (i) and (ii), we choose 
$x:=(x_{s+1},\ldots ,x_{d})$ with 
$x_{i}\in A_{1}$ so that $(\g ,x)$ is a complete system of parameters. 
Let us consider 
the double complex $K^{\bullet}(x^{t}\, ;K^{\bullet}(\g\, ;M))$ whose 
corresponding total complex is $K^{\bullet}(x^{t},\g\, ;M)$. 

Let $\a \in H^{i}_{\im}(M/(\g )M)$ be an element of degree $h$. For
proving (i) (resp. (ii)), we have to show that if $\tau \in 
{\azerc}^{\mu}_{M}$ 
(resp.  $\tau \in {\azer}^{\mu}_{M}$) is such that $\deg \tau >\mu
-i+\s -h$ then  $\tau\a =0\in  H^{i}_{\im}(M/(\g )M)$ for
$\varepsilon\leq i<d-s$ (resp. for $\varepsilon\leq i\leq d-s$). 

For $t:=\max_{i,j}\{ a_{i}(H^{j}(\g\, ;M))\}+\s +s-h+1$ and $\ell\geq
h-\s -s$, we consider the spectral sequence
$$
({'E}_{2}^{pq})_{\geq \ell}=\left[ H^{p}(x^{t}\, ;H^{q}(\g\, ;M))\right] 
_{\geq \ell}\Rightarrow 
H^{p+q}(x^{t},\g\, ;M)_{\geq \ell}.
$$

Now $t>\max_{i}\{ a_{i}(H^{q}(\g\, ;M))\}-\ell$ so that Lemma 2.2 and the
isomorphism  $H^{q}(\g\, ;M)\simeq  H_{s-q}(\g\, ;M)[\s +s]$ provides a
sequence of isomorphisms, 
$$
\left[ H^{p}(x^{t}\, ;H^{q}(\g\, ;M))\right] _{\geq \ell}
\buildrel{\psi_{t,\ell}^{p,q}}\over{\lra} H^{p}_{\im}(H^{q}(\g\,
;M))_{\geq \ell}
\buildrel{\theta_{\ell}^{p,q}}\over{\lra} H^{p}_{\im}(H_{s-q}(\g\,
;M))_{\geq \ell +\s +s}.
$$
As  $H^{p}_{\im}(H_{s-q}(\g\, ;M))=0$ for $p>\varepsilon$  and $q\not=
s$ (notice that $\dim H_{j}(\g\, ;M)\leq \varepsilon$ for $j>0$), we
obtain that $({'E}_{2}^{pq})_{\geq \ell}=0$ for these values of 
$p$ and $q$, so that for $i\geq \varepsilon$ we also have an isomorphism  
$$
 H^{i+s}(x^{t},\g\, ;M)_{\geq \ell}
\buildrel{\phi_{t,\ell}^{i}}\over{\lra}
\left[ H^{i}(x^{t}\, ;H^{s}(\g\, ;M))\right] _{\geq \ell}.
$$

Now for  $i\geq \varepsilon$ the isomorphisms
$\eta^{i,s}_{t,\ell}:=\theta^{i,s}_{\ell}\circ \psi_{t,\ell}^{i,s}\circ
\phi^{i}_{t,\ell}$ provides a commutative diagram  
$$
\xymatrix{ 
H^{i+s}(x^{t},\g\, ;M)_{\geq h-\s -s}
\ar[d]^{\times \tau}\ar[rr]^{\eta_{t,h-\s -s}^{i,s}}&&
[H^{i}_{\im}(M/(\g )M)]_{\geq h}
\ar[d]^{\times \tau}\\
H^{i+s}(x^{t},\g\, ;M)_{\geq \nu -\s -s}
\ar[rr]^{\eta_{t,\nu -\s -s}^{i,s}}&&
[H^{i}_{\im}(M/(\g )M)]_{\geq \nu}\\}
$$
where $\nu :=h+\deg \tau \geq h$, so that $\nu -\s -s=\deg \tau
+h-\s -s>\mu -(i+s)$. 

By hypothesis $[\tau H^{i+s}(x^{t},\g\, ;M)]_{> \mu -(i+s)}=0$
for $i\not= d-s$ (resp. for any $i$), so the
vertical maps are 0, and in particular $\tau \a =0$. 
\fini\medskip

{\bf Proposition 2.3.} {\sl Let $M$ be a finitely generated graded
$R$-module, $x\in R$ an homogeneous element. Assume that
$[xH^{i}_{\im}(M)]_{\mu}=0$ and $\dim (H_{1}(x\, ;M))\leq i$. Then 
$$
H^{i}_{\im}(M)_{\mu}\subseteq H^{i}_{\im}(M/xM)_{\mu}.
$$
}
{\bf Proof.} Set $K:=H_{1}(x\, ;M)=\ker (M[-\de ]\buildrel{\times
x}\over{\lra}M))$ with $\de :=\deg x$. One has an exact sequence 
$$
0\ra K \ra M[-\de ]\buildrel{\times x}\over{\lra} M\ra M/xM\ra 0.
$$
If $[xH^{i}_{\im}(M)]_{\mu}=0$ it gives rise to a surjection,
$$
\xymatrix{
\ker
(H^{i+1}_{\im}(K)_{\mu}\buildrel{can}\over{\lra}H^{i+1}_{\im}(M)_{\mu -\de})
\ar@{->>}[r]&\ker
(H^{i}_{\im}(M)_{\mu}\buildrel{can}\over{\lra}H^{i}_{\im}(M/xM)_{\mu})
\\}
$$
which shows our claim.\fini\medskip

{\bf Proposition 2.4.} {\sl  Let $M$ be a finitely generated graded
$R$-module and $\g =(\g_{1},\ldots ,\g_{t})$
be an homogeneous system of parameters that is a regular sequence on
$M$ outside $V(J)$ with $\dim M/JM\leq 
\varepsilon$.

Let $0\leq \ell <t$, set  $\s_{\ell }:=\deg \g_{1}+\cdots +\deg \g_{\ell }-\ell $ and
assume either that\smallskip 
{\rm (i)}  $\varepsilon\leq i< d-\ell $ and $\g_{j}\in {\azerc}^{\mu}_{M}$
for $j\leq \ell +1$; or \smallskip
{\rm (ii)}  $i\geq \varepsilon$ and $\g_{j}\in {\azer}^{\mu}_{M}$
for $j\leq \ell +1$. \smallskip
\noindent Then,
$$
H^{i}_{\im}(M/(\g_{1},\ldots ,\g_{\ell }))_{\nu}\subseteq H^{i}_{\im}(M/(\g_{1},\ldots ,
\g_{\ell +1}))_{\nu},\quad \forall \nu > \mu -i +\s_{\ell }.
$$

In particular, if  $\g_{j}\in {\azerc}^{\mu}_{M}$
for all $j$ and $\varepsilon\leq i\leq d-t$,
$$
H^{i}_{\im}(M)_{\nu}\subseteq H^{i}_{\im}(M/(\g_{1},\ldots ,
\g_{t}))_{\nu} ,\quad \forall \nu > \mu -i +\s_{t-1}.
$$
}\medskip

{\bf Proof.} For (i), 2.1 (i) implies that $[\g_{\ell +1}H^{i}_{\im}(M/(\g_{1},\ldots ,
\g_{\ell }))]_{\nu}=0$ for $\nu > \mu -i +\s_{\ell }$ so the result follows
from 2.3. The proof of (ii) is similar, replacing 2.1 (i) by  2.1
(ii).\fini\bigskip

{\bf \S 3. Uniform annihilators of homologies}\medskip

{\bf Notation.} If $B$ is a homogeneous quotient of $R$ and $M$  a
finitely generated graded $R$-module, \smallskip
$\bullet$  $J_{B}$  is the jacobian ideal of $B$,\smallskip
$\bullet$  ${\cal T}_{B}$ is the test ideal of $B$ ({\it see} [HH3, \S
3] for details on the characteristic zero case),\smallskip
$\bullet$  ${\ac}_{M}^{*}:=\prod_{i\not= d}{\ac}^{i,-\infty}_{M}$.\medskip

{\bf Proposition 3.1.} {\sl  If $B$ is a homogeneous quotient of
$R$ and $M$  a
finitely generated graded $R$-module,\smallskip 

{\rm (i)} ${\cal T}_{B}\subseteq {\azerc}_{B}^{-\infty}$,\smallskip

{\rm (ii)} ${\ac}^{*}_{M}\subseteq {\azerc}_{M}^{-\infty}
\subseteq {\acc}_{M}^{-\infty}$,\smallskip 

{\rm (iii)} If further $B$ is geometrically reduced and equidimensional,
$J_{B}\subseteq {\cal T}_{B}$.\medskip 

}

{\bf Proof.} The inclusion  ${\cal T}_{B}\subseteq {\azerc}_{B}^{-\infty}$ 
directly follows from the phantom acyclicity criterion [HH1, 9.8].  
The inclusion ${\frak
  C}^{*}_{M}\subseteq {\azerc}_{M}^{-\infty}$ is now standard (see 
[BH, 8.1.3 a]) and ${\azerc}_{M}^{-\infty}\subseteq{\acc}^{i,
-\infty}_{M}$ for $i\not= d$ by 2.1. Finaly, (iii) is a special case of 
[HH3, 3.4].\fini\medskip 

{\bf Proposition 3.2.} {\sl  Let $A$ be a Gorenstein homogeneous 
ring and $I\subseteq A$ an ideal generated by forms of 
degrees $d_{1}\geq\cdots\geq d_{s}$. Assume that $J$ is an 
intersection of primary components of $I$ of codimension $r:=\codim 
(I)$.  Set $\Sc :=\proj (A/J)$ and $\Z :=\proj (A/I)$.

Then there exists an homogeneous 
ideal $\ib$ of $A$, containing $I$, generated in degree at most $\s :=\reg A+d_{1}
+\cdots +d_{r}-r$, such that $\breg (A/J)\leq \s$ and 
$V(\ib )\subseteq W\cup (\Z -\Sc )$, where $W$ is the closed locus of
points of $\Sc$ 
where $\Z$ is not locally a complete intersection F-rational 
scheme.}\medskip 
 
{\bf Proof.} For each point $x\in \Sc$ where $\Z$ is locally a complete
intersection there exists forms $f_{1}^{x},\ldots ,f_{r}^{x}$ of
degrees $d_{1},\ldots ,d_{r}$ that defines a complete intersection
$\Sc \cup \Sc^{x}$ such that $x\not\in \Sc^{x}$. By [CU, 1.7 (iii)],
if $\Sc$ is F-rational at $x$ there exists a form $h_{f}^{x}$ of 
degree $\s$
such that $h_{f}^{x}(x)\not= 0$ and $I_{\Sc}=(f_{1}^{x},\ldots ,f_{r}^{x})
:(h_{f}^{x})$. By [CP, Prop. 2], $A/I_{\Sc}$ is $(\s ,(h^{x}))$-regular
because $A/(f_{1}^{x},\ldots ,f_{r}^{x})$ is $\s$-regular. It
follows that  $A/I_{\Sc}$ is $(\s ,\ib )$-regular with $\ib :=I+\sum_{x,f^{x}}
(f_{1}^{x},\ldots ,f_{r}^{x},h_{f}^{x})$, where the sum is taken over 
$x$ and $f^{x}:=(f_{1}^{x},\ldots 
,f_{r}^{x})$ as above. Notice that $\ib$ is generated in degree at
most $\s$ unless $\s <d_{1}$, in which case $d_{2}=1$ and
$A$ is a polynomial ring, so that choosing $\ib :=A$ gives the result
in this exceptional case. The zero set of $\ib$ is contained by construction  
in the complement of the points of $\Sc$ where  $\Z$ is locally a complete
intersection and  $\Sc$ (or equivalently  $\Z$) is F-rational.\fini\bigskip 

{\bf \S 4. Castelnuovo-Mumford regularity}\medskip

First we note that the Theorem 1 of [CP] readily extends to standard
graded Cohen-Macaulay algebras,\medskip

{\bf Theorem 4.0.}  {\sl  Let $k$ be a field and $A$ a standard graded
Cohen-Macaulay $k$-algebra of dimension $n+1$. If $f_{1},\ldots ,f_{t}$
are forms in $A$ of degrees $d_{1}\geq\cdots\geq d_{t}$ and $\Sc$ is a
subscheme of the zero dimensional component of $\proj (A/(f_{1},\ldots
,f_{t}))$, then
$$
\reg (\Sc )\leq \reg (A)+d_{1}+\cdots +d_{n}-n.
$$
}
{\bf Proof.} See [CP], the proof of Theorem 1.\fini\medskip

We will now extend this result to the components of a scheme of
positive dimension. We have three main results in this
direction.\medskip

{\bf Proposition 4.1.} {\sl Let $k$ be a field and $A$ a standard graded
Gorenstein $k$-algebra, $f_{1},\ldots ,f_{t}$ be forms in $A$ of
degrees $d_{1}\geq\cdots\geq d_{t}\geq 2$, $I$ the ideal they
generate and  $\Z :=\proj (A/I)$. Let $J$ be an intersection of
primary components of $I$ of codimension $r:=\codim (I)$, $\Sc :=\proj 
(A/J)$ and assume that\smallskip

{\rm (1)} $\Z$ is a complete intersection F-rational scheme locally on
$\Sc$ outside finitely many points,\smallskip 
{\rm (2)} $\Sc \cap \overline{\Z -\Sc }=\emptyset$.\smallskip

Then,
$$
\reg (\Sc )\leq (\dim \Sc +1)(\reg (A)+d_{1}+\cdots d_{r}-r-1)+1.
$$
}
\medskip
{\bf Proof.} Consider $\ib$ as given by 3.2. As $\dim (A/(J+\ib ))\leq
1$ by (1), we may choose elements $\g_{1},\ldots ,\g_{t}\in
\ib\subseteq \azer_{A/J}^{\s}$ (with $\s =\reg A+d_{1}+\cdots
+d_{r}-r$) that
forms a sequence of parameters in $A/J$ with $t:=\dim A/J-1$ and $\deg
\g_{i}=\s$. 

By 2.4, $H^{i}_{\im}(A/J)_{\mu}\subseteq
H^{i}_{\im}(A/J+(\g ))_{\mu}$ for $\mu > \s +(t-1)(\s -1)$ and $i\geq
1$. Now 
$H^{i}_{\im}(A/J+(\g ))=0$ for $i\geq 2$ and (2) implies that $\proj
(A/J+(\g ))$ is an isolated component of 
$\proj (A/I+(\g ))$ so that by 4.0, $H^{1}_{\im}(A/J+(\g ))_{\mu}=0$
for $\mu >\s +t(\s -1)$. The claim follows.\fini
\medskip

{\bf Theorem 4.2.} {\sl Let $k$ be a field and $A$ a
  standard graded Gorenstein $k$-algebra, $f_{1},\ldots ,f_{t}$ be
  forms in $A$ of 
degrees $d_{1}\geq\cdots\geq d_{t}\geq 2$, $I$ the ideal they
generate and  $\Z :=\proj (A/I)$. 

Let $J$ be an intersection of
primary components of $I$ of codimension $r:=\codim (I)$ and $\Sc :=\proj 
(A/J)$. Assume that $\dim \Sc$ is positive and\smallskip
{\rm (1)} $\Sc$ is Cohen-Macaulay,\smallskip
{\rm (2)} $\Z$ is a complete intersection locally on $\Sc$ outside
finitely many points,\smallskip 
{\rm (3)} $\Sc$ is F-rational outside a scheme of dimension at most one,\smallskip
{\rm (4)} $\dim (\Z -\Sc )\leq 0$ or $\dim (\Z -\Sc )=1$ and the one
dimensional component of $\Z -\Sc$ is geometrically reduced.\smallskip
Then,
$$
\reg (\Sc )\leq \dim \Sc(\reg (A)+d_{1}+\cdots d_{r}-r-1)+1.
$$
}
\medskip
{\bf Proof.} Consider $\ib$ as given by 3.2. As $\ib \supseteq I$ and
$\dim (A/\ib )\leq 
2$ by (2), (3) and (4), we may choose elements $\g_{1},\ldots 
,\g_{t}\in \ib\subseteq \azer_{A/J}^{\s}$ with $t:=\dim A/J-2$ and
$\deg \g_{i}=\s$, that forms 
a sequence of parameters in $A/J$ and that are such that $\dim
(A/(I+(\g ) ))\leq 2$ and the residual of $\Sc ':=\proj (A/J+(\g ))$ in the top 
dimensional component of $\Z ':=\proj (A/I+(\g ))$ is geometrically reduced. 

By 2.4, $H^{1}_{\im}(A/J)_{\mu}\subseteq   
H^{1}_{\im}(A/J+(\g ))_{\mu}$ for $\mu > \s +(t-1)(\s -1)$. Now
(1) and (4) implies that $\Sc '$ is a union of isolated components
of $\Z'$ of maximal dimension whose residual is reduced so that by
([C, Theorem 35] or [CU, 4.7 (a)(ii)]) $H^{1}_{\im}((A/J+(\g 
))_{\mu}=0$ for $\mu >\s +t(\s -1)$. The claim follows.\fini 
\medskip

The previous propositions only concerns the top dimensional component
of the scheme defined by the given equations. We will now extend our
results to the other isolated components, this time in the
polynomial ring $R$.\medskip
 
{\bf Remark 4.3.} Let $I\subset R$ be a homogeneous ideal, $J$ an
intersection of isolated primary components of $I$ of codimension
$r$. Set $\Z  :=\proj (R/I)$ and $\Sc :=\proj (R/J )$. The following
conditions are equivalent,\smallskip

(i) $\Z$ is smooth locally in codimension $c$ on $\Sc$,\smallskip

(ii) $\Sc$ is smooth outside a scheme of codimension $c+1$, and  
$\ip\in {\rm Ass}(R/I)-{\rm Ass}(R/J)
\Rightarrow \codim(R/(\ip +J))>c$,\smallskip

(iii) $\dim (R/I_{r}(\jac_{R}(I))+J)>c$.\medskip

{\bf Theorem 4.4.} {\sl Let $I\subset R$ be a homogeneous $R$-ideal 
generated by forms of degrees $d_{1}\geq\cdots \geq d_{s}$. Consider
$B:=R/J$, where $J$ is an intersection of isolated primary component 
of codimension $r$ of $I$. Set $\Sc :=\proj (B)$ and assume that,

{\rm (1)}  $\Z :=\proj (R/I)$ have at most isolated singularities on
$\Sc$.\smallskip 

{\rm (2)} $\Sc$ do not meet the other components of $\Z$.
\smallskip
Then
$$
\reg (\Sc )\leq (\dim \Sc +1)(d_{1}+\cdots +d_{r}-r-1)+1.
$$
}

{\bf Proof.} By 4.0 the proposition is true if $\dim \Sc
=0$. Set $B:=R/J$ and assume that $d=\dim B=\dim \Sc +1\geq 2$. Note
that by Remark 4.3,  (1) implies  that there exists $\g_{1},\ldots
,\g_{d-1}\in I_{r}(\jac_{R}(I))\subseteq J_{B}$ that forms a system of
parameters. As $I_{r}(\jac_{R}(I))$ is generated in degree at most $\s
:=d_{1}+\cdots +d_{r}-r$, we may choose $\g_{i}$ such that $\deg
\g_{i}\leq \s$. By 3.1 (i) and (iii), $\g_{i}\in
{\azerc}_{B}^{-\infty}$. The $\g_{i}$'s forms a regular sequence on
the smooth locus of $B$, so that by 2.4 
$$
H^{j}_{\im}(B)\subseteq H^{j}_{\im}(B/(\g_{1},\ldots ,\g_{d-i}))
$$
for $1\leq i<d$ and $1\leq j\leq i$. 

Condition (2) implies that the unmixed part of $\Z _{i}:=\proj
(B/(\g_{1},\ldots ,\g_{d-i}))$ is an isolated component of $\proj
(R/(I+(\g_{1},\ldots ,\g_{d-i})))$, therefore the $a$-invariant of $\Z
_{i}$ is at most $(d-i)\s +d_{1}+\cdots +d_{r}-n-1$ by [CP, corollaire 2]
(note that $\s \geq d_{1}$, unless $d_{2}=1$ where $R/I$ is
Cohen-Macaulay and the result obvious). Therefore
$$
a_{i}(B)\leq a(B/(\g_{1},\ldots ,\g_{d-i}))\leq (d-i)\s +(\s +r)-n-1.
$$
So that $a_{i}(B)+i\leq (d-i+1)\s +i-d$, and
$$
\reg (\Sc )=\max_{i>0}\{ a_{i}(B)+i\} \leq d(\s -1)+1.
$$
\fini\medskip

{\bf Corollary 4.5.} {\sl Let $\Sc \subseteq {\bf P}_{n}$ be a
  projective equidimensional scheme of dimension $d$ over a
  field. Assume that $\Sc$ have at most isolated non smooth points and
  is of  local embedding dimension at most $e$ and set $\kappa :=\min
  \{ n-d,\max \{ e-d,d+1\} \}$. Then,  
$$
\reg (\Sc )\leq \kappa (d+1)(\deg \Sc -1).
$$
If in addition $\Sc$ is Cohen-Macaulay,
$$
\reg (\Sc )\leq \kappa d (\deg \Sc -1).
$$
In particular, if $\Sc$ is smooth, $\kappa =\min \{ n-d,d+1\}$ and 
$$
\reg (\Sc )\leq \kappa d (\deg \Sc -1).
$$
}

Note that a better result for smooth schemes derives from Kodaira
vanishing in characteristic zero ([BM], [BEL] or [CU]).
Also note that the (global) embedding dimension of $\Sc$ is at most 
$\deg \Sc+d-1$, so that a scheme with isolated non smooth points
satisfies $\reg (\Sc )\leq (\dim \Sc +1)(\deg \Sc -1)^{2}$.\medskip

{\bf Proof.} As $\Sc$ have isolated singularities of embedding
dimension  at most $e$, a general linear projection to ${\bf
P}_{m}(k)$ is an isomorphism for $m= \kappa +d$. It
follows that $\reg 
(\Sc )\leq \reg (\pi (\Sc ))$ (and the only possibility for the
inequality to be strict comes from the loss of global sections,
reflected in $H^{1}_{\im}$). Now $\pi (\Sc )$ is scheme defined, up to
points, by equations of degrees at most $\deg \pi (\Sc )=\deg \Sc$
(namely by the cones over $\pi (\Sc )$), and the bound follows from
4.1 or 4.4 for the first result, and from 4.2 for the Cohen-Macaulay
case.\fini\bigskip

\bigskip

{\bf Acknowledgements.} I am very grateful to Patrice Philippon and
Bernd Ulrich for many stimulating discussions. 

\vfill\eject

{\bf References}\medskip

[BH] W. Bruns, J. Herzog, {\it 
Cohen-Macaulay rings}. 
Cambridge Stud. in Adv. Math., 39. 
Cambridge Univ. Press, Cambridge, 1993.\medskip 

[BM] D. Bayer, D. Mumford, {\sl 
What can be computed in algebraic geometry?} 
{\it Computational algebraic geometry and commutative algebra
  (Cortona, 1991), 1--48}. Sympos. Math. XXXIV, Cambridge
Univ. Press, Cambridge, 1993.\medskip  

[BS] D. Bayer, M. Stillman, {\sl  
A criterion for detecting $m$-regularity}, Invent. Math. 
{\bf 87} (1987),  1--11.\medskip 

[BEL] A. Bertram, L. Ein, R. Lazarsfeld, {\sl
Vanishing theorems, a theorem of Severi, and the equations defining
projective varieties}, J. Amer. Math. Soc. {\bf 4} (1991), 587--602.\medskip   

[C]  M. Chardin, {\sl
Applications of some properties of the canonical module in
computational projective algebraic geometry}, J. Symbolic Comput. {\bf
29} (2000), 527--544.\medskip 

[CP] M. Chardin, P. Philippon, {\sl R{\'e}gularit{\'e} et interpolation},
J. Algebraic Geom. {\bf 8} (1999), 471--481. \par
{\it See also the {\rm erratum} at the address} 
http://www.math.jussieu.fr/$\tilde{\;}$chardin/textes.html.\medskip

[CU] M. Chardin, B. Ulrich, {\sl Liaison and Castelnuovo-Mumford regularity}, 
Amer. J. Math., {\it to appear}.\medskip

[HH1]  M. Hochster, C. Huneke, {\sl
Tight closure, invariant theory, and the Brian{\c c}on-Skoda theorem},
J. Amer. Math. Soc. {\bf 3} (1990), 31--116.\medskip
 
[HH2]  M. Hochster, C. Huneke, {\sl
Tight closure in equal characteristic zero}, {\it preprint} (1999).
\medskip

[HH3]  M. Hochster, C. Huneke, {\sl
Comparison of symbolic and ordinary powers of ideals}, 
Invent. Math. {\bf 147} (2002), 349--369.\medskip

[M] C. Miyazaki,  {\sl Graded Buchsbaum algebras and Segre products},
Tokyo J. Math. {\bf 12} (1989),  1--20. \medskip

[MV] C. Miyazaki, W. Vogel, {\sl
Bounds on cohomology and Castelnuovo-Mumford regularity},
J. Algebra {\bf 185} (1996), 626--642.\medskip

[NS1] U. Nagel, P. Schenzel, {\sl Cohomological annihilators and 
Castelnuovo-Mumford regularity},  Contemp. Math. {\bf 159} (1994),
307--328.\medskip

[NS2] U. Nagel, P. Schenzel, {\sl Degree bounds for generators of 
cohomology modules and Castelnuovo-Mumford regularity},  Nagoya
Math. J. {\bf 152} (1998), 153--174.\medskip

[V] W. Vasconcelos, {\it Computational methods in commutative algebra and
algebraic geometry}.  Algorithms and Computation in Mathematics, 2.
Springer-Verlag, Berlin, 1998.\medskip 

\bigskip\bigskip\bigskip

\noindent Marc Chardin, Institut de Math{\'e}matiques,
CNRS \&\ Universit{\'e} Paris 6,

4, place Jussieu, F--75252 Paris {\nrm cedex} 05, France

chardin@math.jussieu.fr

\end